\documentclass[11pt,psfig]{article}
\usepackage{amssymb,amsmath,amscd}
\newtheorem{lem}{Lemma}[section]
\newtheorem{conj}{Conjecture}
\newtheorem{thm}{Theorem}[section]
\newtheorem{cor}{Corallary}[section]
\newtheorem{defn}{Definition}[section]

\newcommand{\f}[1]{\mathfrak{#1}}

\newcommand{\p}{\prime}
\newcommand{\mb}{\mathbb}
\newcommand{\commentout}[1]{}

\newcommand{\mc}{\mathcal}

\newcommand{\Hom}{Hom\,}

\begin{document}
\title{Unitary Representations and Theta Correspondence
for Type I Classical Groups}
\author{Hongyu He  \\
Department of Mathematics \& Statistics \\
Georgia State University \\
email: matjnl@livingstone.cs.gsu.edu \\}
\date{}
\maketitle
\abstract{In this paper, we discuss the positivity of the Hermitian form $(,)_{\pi}$ introduced by Jian-Shu Li in ~\cite{li2}.
Let $(G_1,G_2)$ be a type I dual pair with $G_1$ the smaller group.
Let $\pi$ be an irreducible unitary representation in the semistable range of $\theta(MG_1, MG_2)$ (see ~\cite{theta}). We prove that the invariant Hermitian form $(,)_{\pi}$ is positive
 semidefinite under certain restrictions on the size of $G_2$ and a mild growth condition
on the matrix coefficients of $\pi$. Therefore, if $(,)_{\pi}$ does not vanish, 
$\theta(MG_1, MG_2)(\pi)$ is unitary. \\
\\
Theta correspondence over $\mb R$ was 
established by Howe in (~\cite{howe}). Li showed that theta correspondence 
preserves unitarity for dual pairs in {\it stable range}. 
 Our results generalize the results of Li for type I classical groups (~\cite{li2}). The main result in this paper can be used to construct irreducible unitary representations of classical groups of type I.} 
\section{Introduction}
Let $(G_1,G_2)$ be an irreducible reductive dual pair of type I in $Sp$ (see ~\cite{howe1} ~\cite{li2}).
The dual pairs in this paper will be considered as ordered. For example, the pair $(O(p,q), Sp_{2n}(\mb R))$ is considered different from the pair $(Sp_{2n}(\mb R), O(p,q))$. We will in general assume that the size of $G_1(V_1)$ is less or equal to the size of $G_2(V_2)$. In other words, 
$\dim_D(V_1) \leq \dim_D(V_2)$. Let $Mp$ be the unique double covering of $Sp$. Let $\{1, \epsilon\}$ be the preimage of the identity element in $Sp$.
 For a subgroup $H$ of $Sp$, let
$MH$ be the preimage of $H$ under the double covering. Whenever we use the notation $MH$, $H$ is considered as a subgroup of certain $Sp$. Let $\omega(MG_1, MG_2)$ be a Schr\"odinger model of the oscillator representation of $Mp$. The Harish-Chandra module of $\omega(MG_1, MG_2)$ consists of polynomials multiplied by the Gaussian function. \\
\\
Since the pair $(G_1,G_2)$ is ordered, we use $\theta(MG_1, MG_2)$ to denote the
theta correspondence from $\mc R(MG_1, \omega(MG_1, MG_2))$ to
$\mc R(MG_2, \omega(MG_1, MG_2)$ (see ~\cite{howe}).
In this paper, whenever we talk about "$K$-finite matrix coefficients" or
 "$K$-finite vectors" of a
representation $\pi$ of a real reductive group $G$, "K" is used as a 
generic term for a specified
maximal compact subgroup of $G$. Throughout this paper, we will mainly work within the category of Harish-Chandra modules. A representation of a real reductive group refers to an admissible representation unless stated otherwise. Throughout this paper, a vector in
an admissible representation $\pi$ means that $v$ is in the Harish-Chandra module of $\pi$ which shall be evident within the context.\\
\\
Let $V$ be a vector space of finite dimension. Let $W$ be a subspace of $V$.
A direct complement of $W $ in $V$ is a subspace $U$ such that
$$U \oplus W =V.$$
Now suppose $V$ is equipped with a nondegenerate sesquilinear form $(,)$. The orthogonal complement of $W$ in $V$ consists of 
$$\{ v \in V \mid (v, w)=0 \ \  \forall \ \ w \in W \}.$$
It is denoted by $W^{\perp}$. \\
\\
Let $\pi$ be an irreducible admissible representation of $MG_1$ such that $\pi(\epsilon)=-1$. $\pi$ is said to be in the
semistable range of $\theta(MG_1, MG_2)$ if the function
$$(\omega(MG_1, MG_2) (\tilde g_1) \phi, \psi)(u, \pi(\tilde g_1) v) \qquad (\forall \ \phi, \psi \in 
\omega(MG_1, MG_2); \forall \  u, v \in \pi) $$
is in $L^{1-\delta}(MG_1)$ for all sufficiently small nonnegative $\delta$ (i.e., $ \delta \in [0,c] $ for some $c >0$). We denote the semistable range by $\mc R_s(MG_1, \omega(MG_1, MG_2))$.  Suppose from now on that $\pi$ is in $\mc R_s(MG_1, \omega(MG_1, MG_2))$. For each $\phi, \psi \in \omega(MG_1, MG_2)$  and $u, v \in \pi$, we define an averaging integral 
$$\int_{MG_1}(\omega(MG_1, MG_2) (\tilde g_1) \phi, \psi)(u, \pi(\tilde g_1) v) d \tilde g_1$$
 and denote it by $(\phi \otimes v, \psi \otimes u)_{\pi}$. Thus $(,)_{\pi}$ becomes a real bilinear form on $\omega(MG_1, MG_2) \otimes \pi$. Our definition of $(,)_{\pi}$ differs slightly from the original definition of Li in ~\cite{li2}. Let $\tilde g_2 \in MG_2$ act on $\omega(MG_1, MG_2) \otimes \pi$ by $\omega(MG_1, MG_2)(\tilde g_2) \otimes Id$. In ~\cite{theta}, we show that if $(,)_{\pi} \neq 0$ then $(,)_{\pi}$ descends into a sesquilinear
form on the $K$-finite dual representation of $\theta(MG_1, MG_2) (\pi)$. For $\pi$ unitary, $(,)_{\pi}$ is an invariant Hermitian form on $\theta(MG_1, MG_2)(\pi)$. \\
\\
 For $\pi$ unitary, a conjecture of Li says that 
$(,)_{\pi}$ will always be positive semidefinite (see ~\cite{li2}). If Li's conjecture holds and $(,)_{\pi} \neq 0$, then $\theta(MG_1, MG_2)(\pi)$ is unitary. In this paper, we will prove that $(,)_{\pi}$ is positive semidefinite under certain restrictions. This partly confirms the conjecture of Li. The nonvanishing of certain $(,)_{\pi}$
is proved in ~\cite{thesis} and in ~\cite{non}.\\
\\
We adopt the notations from ~\cite{li2}, ~\cite{theta}, ~\cite{pr1}.
Let $(G_1(V_1), G_2(V_2))$ be a dual pair of type I. Suppose $V_2=V_2^0 \oplus V_2^{\p}$ such that
\begin{enumerate}
\item $(,)_2$ restricted onto $V_2^0$ is nondegenerate;
\item $V_2^{\p}=(V_2^0)^{\perp}$;
\item $V_2^0$ is a direct sum of two isotropic subspaces:
$$V_2^0=X_2^0 \oplus Y_2^0.$$
\end{enumerate}
Obviously, $V_2^0$ will always be of even dimension. Let $X^0=\Hom_D(V_1, X_2^0)$. The oscillator representation $\omega(MG_1(V_1),MG_2(V_2^0))$ can be modeled on $L^2(X^0)$. The action of $MG_1$ on $L^2(X^0)$ is equivalent to the regular action of $G_1$ on $L^2(X^0)$ tensoring with a unitary character $\xi$ of $MG_1$. The generic actions of $G_1$ on $X^0$ are
classified abstractly in Theorem ~\ref{isotropy} and Theorem ~\ref{isotropy2}.\\
\\
Later in this paper, the oscillator representation
$\omega(MG_1(V_1), MG_2(V_2^0))$ is denoted as $\omega(M^0G_1, M^0 G_2^0)$ to indicate the fact that $MG_1(V_1)$ in $(MG_1(V_1), MG(V_2))$ might be different 
from $MG_1(V_1)$ in $(MG_1(V_1), MG_2(V_2^0))$. For the same reason, the oscillator representation $\omega(MG_1(V_1), MG_2(V_2^{\p}))$ is denoted by
$\omega(M^{\p} G_1, M^{\p} G_2^{\p})$.
\begin{thm}[Main Theorem]
Let $(G_1,G_2)$ be a dual pair.
Let $\Xi(g)$ be Harish-Chandra's basic spherical function of $G_1$. Suppose
$\pi$ is an irreducible unitary representation of $MG_1$ such that $\pi(\epsilon)=-1$.
Suppose
 \begin{enumerate} 
\item for any $x,y \in G_1$, the function $\Xi(xgy)$ is
integrable on ${G_1}_{\phi}$ for every
 generic $\phi \in Hom_D(V_1, X_2^0)$ (see Definition ~\ref{gen});
\item the tensor product $\pi_0=\omega(M^{\p}G_1, M^{\p}G_2^{\p}) \otimes \pi \otimes \overline{\xi}$, considered as a representation of $G_1$, is weakly contained in $L^{2}(G_1)$ (see ~\cite{wallach}). 
\end{enumerate}
Then $(,)_{\pi}$ is positive semidefinite. If $(,)_{\pi}$ does
not vanish, then $\theta(MG_1, MG_2)(\pi)$ is unitary.
\end{thm}
{\it Remarks:}
\begin{enumerate}
\item $\omega(M^{\p}G_1, M^{\p}G_2^{\p}))$, $\pi$ and $ \overline{\xi}$ are all projective representations of $G_1$. The fact that $\pi_0$ becomes a unitary
representation of $G_1$ is explained in Part II.
\item The first condition roughly requires that
$$\dim_D(X_2^0) > \frac{\dim_D(V_1)}{2}.$$
 The precise statement depends on the
groups involved.  The function $\Xi(g)|_{{G_1}_{\phi}}$ is in $L^1({G_1}_{\phi})$ implies that
$\Xi(xgy)|_{{G_1}_{\phi}}$ is in $L^1({G_1}_{\phi})$ for any $x, y \in G_1$ and vice versa. 
In fact, $\Xi(g)$ is bounded by a multiple of $\Xi(xgy)$ and vice versa. 
Furthermore, for any compact subset $Y$ of $G_1$, there exists a constant $C$,
such that for any $x, y \in Y$,
\begin{equation}~\label{xi}
\Xi(xg y) \leq C \Xi(g) \qquad (g \in G_1).
\end{equation}
One can prove this by studying the compact picture of the basic spherical principle series representation
(see Chapter VII.1 ~\cite{knapp}). Since this remark may have already been in the literature and a proof will incur a new set of notations, we choose not to give the proof.
\item The growth of matrix coefficients of $\omega(MG_1(V_1), MG_2(V_2^{\p}))$ can be
determined easily. Thus
the second condition can be converted into a growth condition on the matrix coefficients of $\pi$ (see Corollary ~\ref{co1}). 
\item The condition 1 and 2 imply that $\pi$ is in  $\mc R_s(MG_1, \omega(MG_1, MG_2))$. Therefore, $(,)_{\pi}$ is an invariant Hermitian form on $\theta(MG_1, MG_2)(\pi)$. The unitarity of $\theta(MG_1, MG_2)(\pi)$ follows since $(,)_{\pi}$ is positive semidefinite. 
\end{enumerate}
This paper is organized as follows. In Part I, we prove some positivity theorems in the sense of Godement ~\cite{go}.  In Part II, we construct the dual pair 
$(G_1, G_2)$ in terms of homomorphisms and study various subgroups and liftings concerning the tensor decomposition
$$\omega(MG_1, MG_2) \cong \omega(M^0 G_1, M^0 G_2^0) \otimes \omega(M^{\p} G_1, M^{\p} G_2^{\p}).$$
This tensor decomposition is termed as the mixed model in ~\cite{li2}. The interpretation of this tensor product is not completely trivial since $MG_1$, $M^0 G_1$ and $M^{\p} G_1$ may be different double coverings of $G_1$. In Part II, we essentially redo part of section 4 in ~\cite{li2} just to be safe.
In Part III, we study $(\omega(M^0 G_1, M^0 G_2^0), L^2(X^0))$ and classify all 
the generic $G_1$-orbits in $X^0$. This enables us to reduce our averaging integral $(\phi \otimes u, \phi \otimes u)_{\pi}$ to an integral on $G_1$-orbits:
$$\int_{\mc O \in G_1 \backslash X^0} \int_{G_1} \int_{x \in \mc O} \phi(g^{-1} x) \overline{\phi(x)}  (u, \pi_0(g) u) d x dg d [\mc O]. $$
We study each generic orbit integral
$$\int_{G_1} \int_{x \in \mc O} \phi(g^{-1} x) \overline{\phi(x)}  (u, \pi_0(g) u) d x dg $$
 in full generality and convert it into an integral on the isotropic group ${G_1}_x$
$$\int_{{G_1}_x} (\pi_0(g) u_0, u_0) d g.$$
Next, we apply the positivity theorem (Theorem ~\ref{pos}) to show that
this integral is nonnegative. Thus $(,)_{\pi}$ is positive semidefinite.
Finally, we take the pair $(O(p,q), Sp_{2n}(\mb R))$ as an example and state
our main theorem in terms of leading exponents of $\pi$. \\
\\
The author wishes to thank Professors Sigurdur Helgason, Roger Howe, Jian-Shu Li, Tomasz Przebinda, Irving Segal, David 
Vogan Jr. and Nolan Wallach
for encouragements, suggestions and helpful discussions. He also likes to thank Marie Hutchinson for helping him read through the paper of Godement (~\cite{go}).

\section{Part I: Positivity Theorems}
Let $G$ be a real reductive Lie group. Let $K$ be a maximal compact 
subgroup 
of $G$. For any unitary representation $(\pi, H)$ of $G$ and any $\sigma \in \hat K$, let
$H_{\sigma}$ be the $K$-isotypic subspace of $H$. Let $S$ be a subset of $\hat K$. We denote
$$\oplus_{\sigma \in S} H_{\sigma}$$
by $H(S)$. 
\subsection{A generic Theorem}
\begin{thm}~\label{pos0}
Let $G$ be a real reductive Lie group. Let $K$ be a maximal compact 
subgroup 
of $G$. Let $\Xi(g)$ be Harish-Chandra's basic spherical function with respect to $K$.
Let $H$ be a closed unimodular Lie subgroup of $G$. Suppose that $\Xi(g)|_H$ is in $L^1(H)$.
Let $ \phi$ be a positive definite function in $L^{2+\epsilon}(G)(S)$ 
for 
some
finite subset $S$ of
$\hat K$ and any $\epsilon >0$.
Then $\int_H \phi(h) d h \geq 0$.
\end{thm}
Here $L^{2+\epsilon}(G)(S)$ is defined with respect to the left regular action of $G$. \\
\\
Proof: By the GNS construction, we construct a unitary 
representation $(\sigma, \mc H)$
such that $\phi(g)=(\sigma(g) \eta, \eta)$ for some cyclic vector $\eta$ 
in 
$\mc H(S)$.
Since $ \phi$ is a positive definite function in
$L^{2+\epsilon}(G)$ for any $\epsilon >0$, by Theorem 1 in~\cite{chh}, $\sigma$ is weakly contained in $L^2(G)$. 
Thus, there exists a sequence of 
convex linear combinations of diagonal matrix coefficients of $L^2(G)(S)$
$$A_i(g) =\sum_{l=1}^{l_i} a_i^{(l)} (L(g) u_i^{(l)}, u_i^{(l)}), \qquad
\sum_{l=1}^{l_i} a_i^{(l)}=1 \qquad (u_i^{(l)} \in L^2(G)(S), a_i^{(l)} \geq 0)$$ 
such that
$$A_i(g) \rightarrow \phi(g)$$
uniformly on compacta. Let $C_c(G)(S)$ be the space of continuous and compactly supported functions in $L^2(G)(S)$. Since $C_c(G)(S)$ is dense in $L^2(G)(S)$, we choose $u_i^{(l)}$ to be in $C_c(G)(S)$. Notice that
$$A_i(e)= \sum_{l=1}^{l_i}a_i^{(l)} \|u_i^{(l)} \|_{L^2}^2 \rightarrow 
\phi(e)=\|\eta\|^2$$
Hence $\{A_i(e)\}_{i=1}^{\infty}$ is a bounded set. Suppose 
$A_i(e) \leq C$.
From Theorem 2 in ~\cite{chh}, 
$$|(L(g) u_{i}^{(l)}, u_i^{(l)})| \leq  \|u_{i}^{(l)} \|^2_{L^2} (\sum_{\sigma \in 
S} d(\sigma))^{\frac{1}{2}}
\Xi(g).$$
It follows that
\begin{equation}
\begin{split}
|A_i(g)|= & |\sum_{l=1}^{l_i} a_i^{(l)} (L(g) u_{i}^{(l)}, u_i^{(l)})| \\
      \leq &  \sum_{l=1}^{l_i}
a_i^{(l)} \|u_{i}^{(l)} \|^2_{L^2} (\sum_{\sigma \in 
S} d(\sigma))^{\frac{1}{2}} \Xi(g) \\
\leq &  C (\sum_{\sigma \in 
S} d(\sigma))^{\frac{1}{2}} \Xi(g).
\end{split}
\end{equation}
We have proved that $\phi(g)$ can be approximated by positive definite functions $A_i(g)$ such that $A_i(g)$
are uniformly bounded by a fixed multiple of $\Xi(g)$.\\
\\
Now consider the restrictions of $\phi(g)$ to $H$.
From (22.2.3) in ~\cite{dieu},
for $(L(g) u_i^{(l)}, u_i^{(l)})$ with $u_i^{(l)}$ a compactly supported continuous function, 
$$\int_H (L(h) u_i^{(l)}, u_i^{(l)}) d h \geq 0.$$
Thus $\int_H A_i(h) d h \geq 0$. But $A_i(g)|_H$ are bounded by a fixed multiple of an integrable function $\Xi(g)|_H$. By the Dominated Convergence Theorem,
$$\int_H \phi(h) d h = \lim_{i \rightarrow \infty}\int_H A_i(h) d h  \geq 0.$$
Q.E.D.

\subsection{First Variation}

 \begin{thm}~\label{pos1}
 Let $G$ be a real reductive Lie group. Let $K$ be a maximal compact 
subgroup 
of $G$.
Let $H$ be a closed unimodular Lie subgroup of $G$. Let $\Xi(g)$ be
the basic spherical function of $G$ of Harish-Chandra. Suppose that $\Xi(g)|_H$ is in $L^1(H)$. Suppose $(\pi, \mc H)$ is an irreducible unitary representation weakly contained in $L^2(G)$
(see ~\cite{chh}). Let $$v=\sum_{i =1}^{k} \int_{M} \phi_i(x) \pi(\gamma_i(x)) u d x$$ where
\begin{itemize}
\item  $u$ is a $K$-finite vector in $\mc H$;
\item $M$ is a smooth manifold;
\item $\phi_i$ is continuous and is supported on a compact set $X_i \subset M$;
\item  $\gamma_i: M \rightarrow G$ is smooth except a codimension
$1$ subset and the closure of $\gamma_i(X_i)$ is compact.
\end{itemize}
Then 
$$\int_{H} (\pi(h)v, v) d h \geq 0.$$
\end{thm}
The basic idea is to control the function $(\pi(g)v, v)$ by a  convergent integral of left and
right translations of $\Xi(g)$. \\
\\
Proof:
From the proof of Theorem ~\ref{pos0}, we have a sequence
of $K$-finite compactly supported continuous positive definite
functions 
$$A_m(g) \rightarrow (\pi(g)u, u)$$
 uniformly on any compact subset and 
$$|A_m(g)|
\leq C \Xi(g).$$
This implies that
$$|A_m(x g y)| \leq C \Xi(x g y).$$
By the compactness of $supp(\phi_i)$ and the unitarity of $\pi$,
$$(\pi(g)v, v)= \sum_{i,j=1}^k \int_{M \times M}\phi_i(x) \overline{\phi_j(y)}
(\pi(g \gamma_i(x))u, \pi(\gamma_j(y)) u) d x d y .$$
Since the closure of $\gamma_i(X_i)$ is compact,
the closure of $\gamma_j(X_j)^{-1} g \gamma_i(X_i)$ is compact for every $g \in G$. By the Ineqaulity ~\ref{xi}, for any $m$,
\begin{equation}
\begin{split}
 & |\sum_{i,j=1}^k \int_{M \times M}\phi_i(x) \overline{\phi_j(y)}
A_m(\gamma_j(y)^{-1}g \gamma_i(x))) d x d y| \\
\leq &
C \sum_{i,j=1}^k \int_{M \times M}|\phi_i(x)| |\phi_j(y)|
\Xi(\gamma_j(y)^{-1}g \gamma_i(x))) d x d y \\
\leq  & C_1 \Xi(g)
\end{split}
\end{equation} 
for some $C_1>0$. 
Furthermore,
$$ \sum_{i,j=1}^k \int_{M \times M}\phi_i(x) \overline{\phi_j(y)}
A_m(\gamma_j(y)^{-1}g \gamma_i(x)) d x d y \rightarrow (\pi(g)v, v) $$
pointwisely as $m \rightarrow \infty$. 
By the Dominated Convergence Thereom,
$$\int_H (\pi(h) v, v) d h = \lim_{m \rightarrow \infty} \int_H \sum_{i,j=1}^k \int_{M \times M} \phi_i(x) \overline{\phi_j(y)}
A_m(\gamma_j(y)^{-1}h  \gamma_i(x))) d x d y d h.$$
But 
$$A_m(g)= \sum_{l=1}^{l_m} a_m^{(l)}(L(g) u_m^{(l)}, u_m^{(l)}).$$ 
For each
$l$,
\begin{equation}
\begin{split}
 & \int_H \sum_{i,j=1}^k \int_{M \times M} \phi_i(x) \overline{\phi_j(y)}
(L(\gamma_j(y)^{-1}h  \gamma_i(x))u_m^{(l)}, u_m^{(l)}) d x d y d h \\
= & \int_H (L(h) [\sum_{i=1}^k \int_M \phi_i(x)L(\gamma_i(x)) u_m^{(l)} d x],
[\sum_{i=1}^k \int_M \phi_i(x)L(\gamma_i(x)) u_m^{(l)} d x]) d h \\
 \geq  & 0
\end{split}
\end{equation}
because $\sum_{i=1}^k \int_M \phi_i(x) L(\gamma_i(x)) u_m^{(l)} d x$ is a continuous and compactly supported function on $G$.
Hence for every $m$, $$\int_H \sum_{i,j=1}^k \int_{M \times M} \phi_i(x) \overline{\phi_j(y)}
A_m(\gamma_j(y)^{-1}h  \gamma_i(x)) d x d y d h \geq 0.$$
It follows that
$$\int_H (\pi(h)v, v) d h \geq 0.$$
Q.E.D.
\subsection{Second Variation}
\begin{thm}~\label{pos}
 Let $G$ be a real reductive Lie group. Let $K$ be a maximal compact 
subgroup 
of $G$.
Let $H$ be a closed unimodular Lie subgroup of $G$. Let $\Xi(g)$ be
the basic spherical function of $G$ of Harish-Chandra. Suppose that $\Xi( g )|_H$ is in $L^1(H)$. Suppose $(\pi, \mc H)$ is an irreducible unitary representation weakly contained in $L^2(G)$
(see ~\cite{chh}, ~\cite{wallach}). Let $$v=\sum_{i =1}^{n} \int_{M} \phi_i(x) \pi(\gamma_i(x)) u_i d x$$ where
\begin{itemize}
\item  $u_i$ are $K$-finite vectors in $\mc H$;
\item $M$ is a smooth manifold;
\item $\phi_i$ is continuous and is supported on a compact subset $X_i \subset M$;
\item  $\gamma_i: M \rightarrow G$  is smooth except a codimension
$1$ subset and the closure of $\gamma_i(X_i)$ is compact.
\end{itemize}
Then 
$$\int_{H} (\pi(h)v, v) d h \geq 0$$
\end{thm}
The only difference from Theorem ~\ref{pos1} is 
$$ v=\sum_{i =1}^{n} \int_{M} \phi_i(x) \pi(\gamma_i(x)) u_i d x$$
instead of
$$v=\sum_{i =1}^{n} \int_{M} \phi_i(x) \pi(\gamma_i(x)) u d x.$$
Proof: Let $V$ be the linear span of 
$$\{ \pi(k) u_i \mid i \in [1,n], k \in K \}.$$
Since $u_i$ are $K$-finite, $V$ is a finite dimensional representation of $K$.
Let $u$ be a $K$-cyclic vector in $V$. Let $C(K)$ be the space of continuous functions on $K$. Consider the action of $C(K)$ on $u$:
$$\pi(f)u = \int_{K} f(k) \pi(k) u d k. $$
 Apparently, $\pi(C(K))u = V$. Let
$$u_i = \int_{K} f_i(k) \pi(k) u d k.$$
Then
\begin{equation}
\begin{split}
 v = & \sum_{i =1}^{n} \int_{M} \phi_i(x) \pi(\gamma_i(x)) u_i d x \\
 = & \sum_{i =1}^{n} \int_{M} \phi_i(x) \pi(\gamma_i(x)) \int_K f_i(k) \pi(k) u d k  d x \\
= & \sum_{i =1}^{n} \int_{M} \int_K \phi_i(x) f_i(k) \pi(\gamma_i(x) k) u d x dk
\end{split}
\end{equation}
Apply Theorem ~\ref{pos1} to functions $\phi_i(x) f_i(k)$ on $M \times K$ and
$$\gamma_i^*: (x, k) \in M \times K \rightarrow \gamma_i(x) k \in G.$$
The conclusion follows immediately. Q.E.D. 
\begin{conj}
Let $G$ be a real reductive group. Let $K$ be a maximal compact subgroup of $G$. Let $\Xi(g)$ be Harish-Chandra's basic spherical function. Let $H$ be a subgroup of $G$ such that $\Xi(g)|_H$ is in $L^1(H)$. Let $\phi(g)$ be a positive definite continuous function bounded by $\Xi(g)$. Then
$\int_H \phi(h) d h \geq 0$.
\end{conj}

\section{Part II: Dual Pairs and Mixed Model}
The basic theory on the mixed model of the oscillator representation is covered in ~\cite{li2} with reference to an unpublished note of Howe. We redo part of section 4 of ~\cite{li2} with emphasis on the actions of various coverings of $G_1$ regarding the mixed model
$$\omega(MG_1, MG_2) \cong \omega(M^0 G_1, M^0 G_2^0) \otimes \omega(M^{\p} G_1, M^{\p} G_2^{\p}).$$
Let $V_1$ be a vector space over $D$ equipped with a sesquilinear form
$(,)_1$, $V_2$ be a vector space over $D$ equipped with a sesquilinear form $(,)_2$. Suppose one sesquilinear form is $\sharp$-Hermitian and the other is $\sharp$-skew Hermitian. Let $G_i$ be the isometry group of $(,)_i$. Let $V= Hom_D(V_1, V_2)$ be the space of $D$-linear homomorphisms from $V_1$ to $V_2$. 

\subsection{Setup}
Let $\phi, \psi \in V$, $v_1, u_1 \in V_1$ and $v_2 \in V_2$. We define a unique
$\phi^*(v_2)$ such that
$$(\phi^*(v_2), v_1)_1= (v_2, \phi(v_1))_2.$$
It is easy to verify that $\phi^* \in Hom_D(V_2,V_1)$.
Thus we obtain a $*$ operation from $V$ to
$V^*=Hom_D(V_2,V_1)$. Let $a \in \mb R$. Then
$$((a \phi)^*(v_2), v_1)_1=(v_2, a \phi(v_1))_2=a(v_2, \phi( v_1))_2=a(\phi^*(v_2), v_1)_1
=(a \phi^*(v_2), v_1)_1.$$
Therefore, the $*$-operation is real linear.\\
\\
 Let $tr(*)$ be the real trace of a real linear endomorphism. Since $V$ and $V^*$ are real vector spaces, we can now define a real bilinear form $\Omega$ on $V$ as follows
$$\Omega(\phi, \psi)=tr( \psi^* \phi).$$
We observe that
$$(\psi^* \phi(v_1), v_1^{\p})_1=(\phi(v_1), \psi(v_1^{\p}))_2
=\pm (\psi(v_1^{\p}), \phi(v_1))^{\sharp}_2= \pm (\phi^* \psi(v_1^{\p}), v_1)^{\sharp}_1 = -(v_1, \phi^* \psi(v_1^{\p}))_1. $$
Define a $*$-operation on $End_D(V_1)$ by
$$(A^* u_1, v_1)_1=(u_1, A(v_1))_1 \qquad (\forall \, \, A \in End_D(V_1)).$$
Then,
$(\phi^* \psi)^*= -\psi^* \phi$.
It follows that
$$\Omega(\psi, \phi)=tr(\phi^* \psi)=tr((\phi^* \psi)^*)=tr(-\psi^* \phi)=-\Omega(\phi, \psi).$$
It is easy to verify that $\Omega$ is nondegenerate. 
Therefore, $\Omega$ is a real symplectic form on $V$. \\
\\
Next we define the action of $G_1$ on $V$ as follows
$$ (g_1 \phi) (v_1)= \phi(g_1^{-1} v_1) .$$
We observe that
\begin{equation}
\begin{split}
 & ((g_1 \psi)^* (g_1 \phi) (u_1), v_1)_1 \\
= & ( (g_1 \phi) (u_1), (g_1 \psi)(v_1))_2 \\
= & (\phi(g_1^{-1} u_1), \psi (g_1^{-1} v_1))_2 \\
= & (\psi^* \phi (g_1^{-1} u_1), g_1^{-1} v_1)_1 \\
= & (g_1 (\psi^* \phi) g_1^{-1} u_1, v_1)_1.
\end{split}
\end{equation}
It follows that 
$$\Omega(g_1 \phi, g_1 \psi)=tr((g_1 \psi)^* (g_1 \phi))=tr (g_1 \psi^* \phi g_1^{-1})=tr(\psi^* \phi)= 
\Omega(\phi, \psi).$$
Therefore $G_1$ is in $Sp(V, \Omega)$.
We define the action of $G_2$ on $V$ similarly by
$$(g_2 \phi)(v_1)= g_2 \phi(v_1).$$
One can verify 
 that $G_2$ also preserves $\Omega$.
In addition, the action of $G_1$ commutes with the action of $G_2$.
\subsection{Subgroups}
Let $V_2^0$ be a $D$-linear subspace of $V_2$ such that
\begin{itemize}
\item 
$(,)_2$ restricted to $V_2^0$ is nondegenerate;
\item 
There exist isotropic subspaces $X_2^0$ and $Y_2^0$ such that
$$X_2^0 \oplus Y_2^0 = V_2^0$$
\end{itemize}
Let $V_2^{\p}$ be the space of vectors perpendicular to $V_2^0$ with
respect to $(,)_2$.
Write 
$$X^0= Hom_D(V_1, X_2^0), \qquad Y^0=Hom_D(V_1, Y_2^0),$$
$$V^{\p}=Hom_D(V_1, V_2^{\p}), \qquad V^0=Hom_D(V_1, V_2^0).$$
For any $\phi, \psi \in X^0$,
$$(\psi^*\phi v_1, u_1)_1=(\phi v_1, \psi u_1)_2=0 \qquad (v_1, u_1 \in V_1).$$
Thus $\Omega(\psi, \phi)=tr(\phi^* \psi)=0$.
$X^0$ is an isotropic subspace of $(V, \Omega)$. For the same reason, $Y^0$ is also an isotropic subspace of $(V, \Omega)$.
Furthermore, we have
$$V=V^{\p} \oplus V^0, \qquad V^0=X^0 \oplus Y^0.$$
Let $G_2^0$ be the subgroup of $G_2$ such that $G_2^0$ restricted to
$V_2^{\p}$ is trivial. Then $G_2^0$ is isomorphic to $G_2(V_2^0)$.
Let $G_2^{\p}$ be the subgroup of $G_2$ such that $G_2^{\p}$ restricted to $V_2^0$ is trivial. Then $G_2^{\p}$ is isomorphic to $G_2(V_2^{\p})$.\\
\\
Let $\Omega^0$ be the restriction of $\Omega$ on $V^0$. Let $\Omega^{\p}$ be the restriction of $\Omega$ on $V^{\p}$.
Then $Sp(V^0, \Omega^0)$ and $Sp(V^{\p}, \Omega^{\p})$ can be embedded into $Sp(V, \Omega)$ diagonally.
Let $GL(X^0, Y^0)$ be the subgroup of $Sp(V^0, \Omega^0)$ stabilizing
$X^0$ and $Y^0$. Since $G_1$ and $G_2^0$ act on  $V^0$,
we obtain a dual pair 
$$(G_1, G_2^0) \subseteq Sp(V^0, \Omega^0).$$
We denote this embedding by $i^0$. On the other hand,
since $G_1$ and $G_2^{\p}$ act on $V^{\p}$, we obtain another dual pair
$$(G_1, G_2^{\p}) \subseteq Sp(V^{\p}, \Omega^{\p}).$$
We denote this embedding by $i^{\p}$.
Now the group $G_1$ is embedded into $Sp(V, \Omega)$ by $i^0 \times i^{\p}$.
We denote this embedding by $i$.

\subsection{Metaplectic Covering and Compatibility}
For any symplectic group $Sp$, there is a unique double covering $MSp$. We call this
the metaplectic covering. Let $\epsilon$ be the nonidentity element in $MSp$ whose image
is the identity element in $Sp$. For any subgroup $G$ of $Sp$, let $MG$ be the preimage of
$G$ under the metaplectic covering. Then every $MG$ contains $\epsilon$. \\
\\
 Let $M^0Sp(V^0, \Omega^0)$,
$M^{\p}Sp(V^{\p}, \Omega^{\p})$ and $MSp(V, \Omega)$ be the metaplectic coverings
of $Sp(V^0, \Omega^0)$,
$Sp(V^{\p}, \Omega^{\p})$ and $Sp(V, \Omega)$ respectively. Let $M^0$, $M^{\p}$
and $M$ be the covering maps respectively. When we consider $Sp(V^0, \Omega^0)$ as a subgroup of $Sp(V, \Omega)$, we obtain a group $MSp(V^0, \Omega^0)$. On the other hand, $Sp(V^0, \Omega^0)$ has its own metaplectic covering, namely, $M^0 Sp(V^0, \Omega^0)$.

\begin{lem}[compatibility]
The group $MSp(V^0, \Omega^0)$ is isomorphic to
$M^0Sp(V^0, \Omega^0)$.
\end{lem}
Proof: It suffices to show that $MSp(V^0, \Omega^0)$ does not split.
Suppose $MSp(V^0, \Omega^0)$ splits.
Let $K$ be a maximal compact subgroup of $Sp(V, \Omega)$ such that $K^0=K \cap Sp(V^0, \Omega^0)$
is a maximal compact subgroup of $Sp(V^0, \Omega^0)$. Then $MK^0$ splits.
On the other hand, $K$ can be identified with a unitary group $U$. The metaplectic 
covering of $U$ can be represented by
$$\{ (\xi, g) \mid \xi^2=\det g, g \in U \}$$
For the subgroup $K^0$, we see that $MK^0$ must be the nontrivial double covering
of $K^0$. It does not split. We reach a contradiction. Q.E.D.\\
\\
This Lemma basically asserts that if a smaller symplectic group is embedded in a bigger symplectic group canonically, then the metaplectic covering on the smaller group is compatible with the metaplectic covering on the bigger group.
Let $$\tilde{i^0}:(M^0 G_1, M^0 G_2^0) \subseteq M^0Sp(V^0, \Omega^0)$$
be the lifting of $i^0$.
Let $$\tilde{i^{\p}}:(M^{\p} G_1, M^{\p} G_2^{\p}) \subseteq M^{\p} Sp(V^{\p}, \Omega^{\p})$$
be the metaplectic
lifting of $i^{\p}$.
Let $$\tilde{i}:(MG_1, MG_2) \subseteq MSp(V, \Omega)$$
be the lifting of $i$.
According to the compatibility lemma, we may consider $M^0 Sp(V^0, \Omega^0)$ 
and $M^{\p}Sp(V^{\p}, \Omega^{\p})$ as subgroups
of $MSp(V, \Omega)$. These two subgroups intersect. The intersection is
$\{ 1, \epsilon \}$. \\
\\
Consider the natural multiplication map
$$j: M^0 Sp(V^0, \Omega^0) \times M^{\p}Sp(V^{\p}, \Omega^{\p}) \rightarrow
MSp(V, \Omega).$$
Its kernel is $\{(1,1), (\epsilon, \epsilon) \}$. If $g \in G_1$, then
$$i(g)=(i^0(g), i^{\p}(g)) \in Sp(V^0, \Omega^0) \times Sp(V^{\p}, \Omega^{\p})
 \subseteq Sp(V, \Omega)$$
The covering group $MG_1$ is then isomorphic to the quotient
$$\{ j(g^0, g^{\p}) \mid g^0 \in M^0G_1, g^{\p} \in M^{\p}G_1, M^0(g^0)=g=M^{\p} (g^{\p}) \}/\{(1,1), (\epsilon, \epsilon) \}. $$

\begin{lem}
Each element in $MG_1$ can be expressed as 
$j(g^0, g^{\p})$ with 
$$(g^0 \in M^0G_1, g^{\p} \in M^{\p}G_1, M^0(g^0)=M^{\p} (g^{\p}))$$
up to a factor of
$$\{(1,1), (\epsilon, \epsilon)\}.$$
\end{lem}

\begin{lem}~\label{tensor0} As a group, 
$$M^0 G_1 \cong \{ (g, g^{\p}) \mid M(g)=M^{\p}(g^{\p}), g \in MG_1, g^{\p} \in M^{\p} G_1 \}/ \{(1,1), (\epsilon, \epsilon) \}.$$
\end{lem}

\subsection{Oscillator Representation as tensor product}
\begin{thm}~\label{tensor}
The representation
$$\omega(M^0G_1, M^0 G_2^0) \otimes \omega(M^{\p} G_1, M^{\p} G_2^{\p})$$
restricted to
$$\{ j(g^0, g^{\p}) \mid g^0 \in M^0G_1, g^{\p} \in M^{\p}G_1, M^0(g^0)=g=M^{\p} (g^{\p}) \} $$
descends into $ \omega(MG_1, MG_2)|_{MG_1}$. 
\end{thm}
Proof: Suppose $g \in MG_1$. Then $g$ can be written as
$$(g^0, g^{\p}) \mid g^0 \in M^0G_1, g^{\p} \in M^{\p}G_1, M^0(g^0)=M^{\p} (g^{\p})$$
up to a multiplication by
$$\{(1,1), (\epsilon, \epsilon)\}.$$
It is easy to see that
$$\omega(MG_1, MG_2)(1,1)=id=\omega(M^0 G_1, M^0 G_2^0)(\epsilon) \otimes \omega(M^{\p}G_1,
M^{\p} G_2^{\p})(\epsilon)$$
It follows that 
$$\omega(MG_1, MG_2)(g)= \omega(M^0 G_1, M^0 G_2^{0})(g^0) 
\otimes \omega(M^{\p} G_1, M^{\p} G_2^{\p}) (g^{\p})$$
Our theorem is proved. Q.E.D. \\
\\
Let $\pi$ be an irreducible unitary representation of $MG_1$ in the semistable
range of $\theta(MG_1, MG_2)$ such that $\pi(\epsilon)=-1$.
Identify the representation $\omega(MG_1, MG_2)^c \otimes \pi$ with
$$\omega(M^0G_1, M^0 G_2^0)^c \otimes (\omega(M^{\p} G_1, M^{\p} G_2^{\p})^c\otimes \pi).$$
From
Lemma ~\ref{tensor0}, $g^0 \in M^0 G_1$ can be represented by a pair $(\tilde g, g^{\p})$ up to a multiplication of $(\epsilon, \epsilon)$. Since
$$\omega(M^{\p} G_1, M^{\p} G_2^{\p})^c(\epsilon) \pi(\epsilon)=id,$$
we can write
$$(\omega(M^{\p} G_1, M^{\p} G_2^{\p})^c\otimes \pi) (g^0)=\omega(M^{\p} G_1, M^{\p} G_2^{\p})^c ( g^{\p}) \otimes \pi(g).$$
The proof of Theorem ~\ref{tensor} shows that
$$\omega(M^{\p} G_1, M^{\p} G_2^{\p})^c \otimes \pi$$
can be regarded as a unitary representation of $M^0 G_1$. 
\subsection{Schr\"odinger Model of $\omega(M^0G_1, M^0 G_2^0)$}
Recall $V^0=X^0 \oplus Y^0$ and both $X^0, Y^0$ are Lagrangian in $(V^0, \Omega^0)$. Let $GL(X^0,Y^0)$ be the subgroup of $Sp(V^0, \Omega^0)$  stabilizing $X^0$ and $Y^0$. Then
$$GL(X^0, Y^0) \cong GL(X^0) \cong GL(Y^0).$$
Let $L^2(X^0)$ be a Schr\"odinger model of $\omega(M^0G_1, M^0 G_2^0)$ (see ~\cite{theta}, ~\cite{wallach0}).
The group $M^0GL(X^0,Y^0)$ acts on $L^2(X^0)$ naturally. Since
$G_1$ is a subgroup of $GL(X^0, Y^0)$,
an element in the group $M^0 G_1$ can be written as 
$$(\xi, g) \mid g \in G_1, \xi \in \mb C $$
such that the operator
$$(\omega(M^0G_1, M^0 G_2^0)(\xi,g) \phi)(x)= {\xi} \phi(g^{-1} x) \qquad (x \in X^0, \phi \in L^2(X^0))$$
is unitary. \\
\\
Consider
\begin{equation}~\label{general0}
\int_{M^0 G_1} (\omega(M^0 G_1, M^0 G_2^0)(\xi,g) \phi, \psi)(u, (\omega(M^{\p}G_1, M^{\p}G_2^{\p})^c \otimes \pi) (\xi, g) v) 
d g d \xi
\end{equation}
with $u, v \in \omega(M^{\p}G_1, M^{\p}G_2^{\p}) \otimes \pi$.
Since the group action of  $G_1$ on $L^2(X^0)$ is already unitary, $\xi$ is
 a unitary character of $M^0 G_1$.
 Thus
$\overline{\xi} \otimes \omega(M^0G_1, M^0G_2^0)$ can be viewed as  a 
unitary representation of $G_1$. 
Moreover, 
$$\overline{\xi} {\omega(M^0G_1, M^0G_2^0)}(g,\xi)\phi(x)=\phi(g^{-1}x).$$
Define
$$\pi_0= \overline{\xi} \otimes (\omega(M^{\p} G_1, M^{\p} G_2^{\p})^c \otimes \pi).$$ 
Viewing $(\omega(M^{\p} G_1, M^{\p} G_2^{\p})^c \otimes \pi)$ as a representation of $M^0 G_1$, $\pi_0$ descends into a unitary representation of $G_1$. \\
\\
Tensor products with $\overline \xi$ here do not change
the ambient spaces. However, the group actions differ by a unitary character.
Now, the integral (~\ref{general0}) becomes a multiple of
\begin{equation}~\label{general}
\int_{G_1} \int_{X^0} \phi(g^{-1} x) \overline{\psi(x)} d x (u, \pi_0(g) v) dg.
\end{equation}
This integral can be expressed as orbital integral 
$$\int_{G_1} \int_{\mc O \in G_1 \backslash X^0} \int_{x \in \mc O}. $$
In Part III, we will classify the generic $G_1$-orbits in $X^0$ and study each generic orbital integral
$$\int_{G_1} \int_{x \in \mc O} \phi(g^{-1} x) \overline{\psi(x)}  (u, \pi_0(g) v) d x d g. $$

\section{Part III: Orbital Integrals}

Recall that $X^0=Hom_D(V_1, X_2^0)$.
We need to classify the orbital
structure of the $G_1$-action on $X^0$. Let
$m= \dim_D V_1 $ and $ \dim_D X_2^0 =p$. If
$m \leq p$, $(G_1, G_2)$ is said to be in the stable range. The action of $G_1$ on $X^0$ is almost free.
This case is already treated
in ~\cite{li2}. For $(G_1, G_2)$ in the stable range, our approach can be simplified and in deed coincides with Li's approach in ~\cite{li2}.
From now on, assume $m \geq p$.
The set of nonsurjective homomorphisms from $V_1$ to $X_2^0$ is of measure zero. 
Hence we will
focus on surjective homomorphisms in $X^0$. We denote the set of surjective homomorphisms by $X_0^0$.
Let $\phi \in X_0^0$.
\subsection{The Isotropic subgroup ${G_1}_{\phi}$}

Let $e_1, e_2, \ldots, e_m$ be a $D$-linear basis for $V_1$, 
and $f_1, f_2, \ldots, f_p$
be a $D$-linear basis for $X_2^0$. Then $\phi$ is uniquely determined by
$$\phi(e_1),\phi(e_2), \ldots, \phi(e_m).$$
We will determine the "generic" isotropic subgroups of the $G_1$-action on $X^0_0$.
Suppose $g \in G_1$ stabilizes $\phi$. In other words,
$$\phi(u) = (g \phi)(u)=\phi(g^{-1} u) \qquad (\forall \ \ u).$$
This implies that $\ker(\phi)$ is stabilized by $g$. Therefore $\ker(\phi)^{\perp}$ is also
stabilized by $g$. 
\begin{lem}Let $g \in G_1$ and $\phi \in X_0^0$. Then
$\phi$ is fixed by $g $ if and only if any vector in $\ker(\phi)^{\perp}$ is 
fixed by $g$. 
\end{lem}
Proof: Suppose $\phi$ is fixed by $g$. Let $(v, \ker \phi)_1=0$. 
We choose an arbitrary $u \in V_1$. Since $\phi(g^{-1} u)=\phi(u)$, 
$g^{-1}u-u \in \ker \phi$.
This implies that
$(v, g^{-1}u-u)_1=0$.
Thus
$(gv, u)_1=(v, u)_1$ for every $u \in V_1$. It follows that $gv=v$. $g$ fixes every vector in $v \in \ker \phi^{\perp}$. \\
\\
Conversely, suppose $gv=v$ for any $(v, \ker \phi)=0$. We choose an arbitrary
 $u \in V_1$. Then 
$(g v-v, u)_1=0$.
 Hence
$(v, g^{-1} u-u)_1=0$
for every $v \in \ker \phi^{\perp}$. From the nondegeneracy of $(,)_1$,
$$ g^{-1} u -u \in (\ker \phi^{\perp})^{\perp}= \ker \phi$$
Therefore, $\phi(g^{-1} u-u)=0$ for every $u \in V_1$. It follows that
$g \phi =\phi$. Q.E.D.

\begin{thm}~\label{isotropy} Let $\phi$ be a surjective homomorphism from
$V_1$ to $X_2^0$. Then the isotropic subgroup ${G_1}_{\phi}$ is the subgroup that fixes all vectors in $\ker(\phi)^{\perp}$. 
\end{thm}
The restriction of $(,)_1$ onto $\ker \phi^{\perp}$ contains a null space, namely, 
\begin{equation}~\label{e1}
W=\ker \phi \cap \ker \phi^{\perp}.
\end{equation}
$W$ is an isotropic subspace of $V_1$ and it may or may not be trivial.
Let $U$ be a direct complement of $W$ in $\ker \phi^{\perp}$, i.e.,
\begin{equation}~\label{e2}
U \oplus W= \ker \phi^{\perp}.
\end{equation}
 Then 
$(,)_1$ restricted to $U$ is nondegenerate. Thus 
$(,)_1$ restricted onto $U^{\perp}$ is a nondegenerate sesquilinear form. Since the group ${G_1}_{\phi}$ fixes all vectors in $\ker \phi^{\perp}$ and $U \subseteq \ker \phi^{\perp}$, ${G_1}_{\phi}$ can be identified with the subgroup of $G_1(U^{\perp}$ that fixes all vectors in $W$. \\
\\
From Equations ~\ref{e1}, ~\ref{e2}, $\ker \phi$ is the orthogonal complement of $W$ in $U^{\perp}$. From Equations (28),(29) in ~\cite{li2}, ${G_1}_{\phi}$ is a twisted product of
$G_1(\ker \phi /W)$ with a at most two-step nilpotent group $N$.
\begin{thm}~\label{isotropy2}
For orthogonal groups, we take $G_1=SO(p,q)$. The isotropic subgroup ${G_1}_{\phi}$ is a twisted product of a classical group of the same type with a at most two-step nilpotent group $N$. It is always unimodular. 
\end{thm}
Proof: To show that ${G_1}_{\phi}$ is unimodular, one must show that the adjoint action of $G_1(\ker \phi /W)$ on
the Lie algebra $\f n$ has determinant $1$. 
This is obvious since $\f n$ as a $G_1$-module decomposes into
direct sum of of trivial representations and 
the standard representations. Q.E.D. 
\subsection{Generic Element}

The homomorphism $\phi$ induces an isomorphism 
$$[\phi]: V_1/\ker \phi \rightarrow X_2^0$$
Notice that $\ker \phi $ can be regarded as a point in the Grassmannian 
$\mc G(m,m-p)$.
We obtain a fibration
$$GL_p(D) \rightarrow X_0^0 \rightarrow \mc G(m,m-p).$$
The projection maps $\phi$ to $ \ker \phi$. The fiber contains all isomorphisms
from $V_1/\ker \phi$ to $X_2^0$. Thus the fiber can be identified with $GL_p(D)$.

\begin{defn}~\label{gen}
{\it Generic elements} in $X^0$ are those surjective $\phi$ such that
\begin{enumerate}
\item either $(,)_1$ restricted on $\ker(\phi)$ is nondegenerate;
\item or if the above case is not possible,
$$\dim_D(\ker(\phi) \cap \ker(\phi)^{\perp})=1$$
\end{enumerate}
Let $X_{00}^0$ be the subset of generic elements. The subspaces $\ker(\phi)$ in
$\mc G(m,m-p)$ are called generic $(m-p)$-subspaces. The set of generic $(m-p)$-subspaces 
is denoted by
$\mc G_0(m,m-p)$.
\end{defn}
Consider the following fibration
$$GL_p(D) \rightarrow X_{00}^0 \rightarrow \mc G_0(m,m-p).$$
Since the set $\mc G_0(m,m-p)$ is open and dense in $\mc G(m,m-p)$, the set $X_{00}^0$
is open and dense in $X_0^0$. Therefore, $X_{00}^0$ is open and dense in $X^0$.\\
\\
First, suppose $(,)_1$ restricted to $\ker(\phi)$ is nondegenerate.
 We must have
$$\ker(\phi) \oplus \ker(\phi)^{\perp}=V_1$$
 The isotropic subgroup ${G_1}_{\phi}$ can be identified with $G_1(\ker(\phi))$ by restriction according to Theorem ~\ref{isotropy}. It 
is a smaller group of type $G_1$. The group ${G_1}_{\phi}$ is automatically unimodular.\\ 
\\
Secondly, suppose 
$$\dim_D(\ker(\phi) \cap \ker(\phi)^{\perp})=1$$
Notice that this case does not occur for $O(p,q)$.
From Theorem ~\ref{isotropy2}, ${G_1}_{\phi}$ is a unimodular group. We obtain
\begin{cor}~\label{isotropy3} 
For any generic element $\phi \in X_{00}^0$, the isotropy subgroup ${G_1}_{\phi}$ is 
always unimodular.
\end{cor}
\subsection{Averaging Integral Revisited}

Let $\pi$ be an irreducible unitary representation in the semistable range of
$\theta(MG_1, MG_2)$. Recall that
 $$ \pi_0 =\omega(M^{\p} G_1, M^{\p} G_2^{\p})^c \otimes \pi \otimes \overline{\xi}$$
is a unitary representation of $G_1$.
Consider the integral
\begin{equation}~\label{average}
\int_{G_1} \int_{X^0} \phi(g^{-1} x) \overline{\psi(x)} d x (u, \pi_0(g) v) dg
\end{equation}
where $\phi, \psi$ are $K$-finite vectors in $L^2(X^0)$ and $u, v \in \pi_0$.
\begin{thm}~\label{conver} Let $\pi$ be an irreducible unitary representation in the 
semistable range of
$\theta(MG_1, MG_2)$. Let $\phi, \psi$ be in the Harish-Chandra module of
$\omega(M^0 G_1, M^0 G_2^0)$. Let $ u, v \in \pi_0$.
 Then the function
$\phi(g^{-1} x) \overline{\psi(x)}  (u, \pi_0(g) v)$
is continuous and absolutely integrable on $G_1 \times X^0$.
Therefore we have
$$\int_{G_1} \int_{X^0} \phi(g^{-1} x) \overline{\psi(x)} d x  (u, \pi_0(g) v) dg
= \int_{X^0} \int_{G_1} (\phi(g^{-1} x) \overline{\psi(x)}  (u, \pi_0(g) v) d g d x$$
\end{thm}
From our discussion in Part II, the integral ~\ref{average} is a form of the averaging integral under the mixed model
 $$\omega(MG_1, MG_2) \cong \omega(M^0 G_1, M^0 G_2^0) \otimes \omega(M^{\p} G_1, M^{\p} G_2^{\p}).$$
The absolute integrabilty of $\phi(g^{-1} x) \overline{\psi(x)}  (u, \pi_0(g) v)$ is guaranteed by the semistable condition (see ~\cite{theta}). We skip the proof.

\subsection{Orbital integral in General}
First, let me quote a simplified version of Theorem 8.36 from ~\cite{knapp1}.
\begin{thm}~\label{int}
Let $G$ be a unimodular group and $H$ be a closed unimodular subgroup of $G$.
Let $d g$ and $d h$ be their Haar measures, respectively.
Then up to a scalar, there exists a unique $G$-invariant measure $d [g H]$ on $G/H$. Furthermore, this measure can be normalized such that for any $L^1$ function on $G$,
$$\int_G f(g) d g =\int_{G/H} \int_H f(g h) d h d [ g H].$$
\end{thm}
Suppose $\tau$ is a unitary representation of $G$, $u$ and $v$ are
$K$ finite vectors in $\tau$. 
\begin{thm}~\label{orbit0} Let $G$ be a  real reductive group, and $M$ be a $G$-homogeneous space. 
\begin{itemize}
\item Let $x_0$ be a fixed base point and $G_0$ be the isotropy group of $x_0$. Suppose that $G_0$ is
unimodular. Then $M$ is isomorphic to $G/G_0$ and possesses a $G$-invariant measure. 
\item Let $\gamma: M \rightarrow G$ be a smooth section of the principle
bundle
$$B: G_0 \rightarrow G \rightarrow M$$
except for a subset of at least codimension $1$. Assume $\phi(y)$ is an absolutely integrable function on $M$. Then
$$v_0=\int_{M} \overline{\phi(y)} \tau(\gamma(y)^{-1}) v d y  $$
is well-defined.
\item
Assume $\phi(g^{-1} x_0)  (u, \tau(g) v)$ is integrable as a function on $G$.
Then we have
$$ \int_{G} \phi(gx_0) (\tau(g) u, v) d g =\int_{G_0}(\tau(g_0) u, v_0) d g_0.$$
\end{itemize}
\end{thm}
Proof: (1) follows directly from Theorem ~\ref{int} by identifying $M$ with $G/G_0$. \\
\\ 
Since $\tau$ is unitary and $\phi(y)$ is integrable, $v_0$ is well-defined. (2) is proved. \\
\\ Notice that 
$\gamma(y) G_0 x_0 =y$.
We compute
\begin{equation}
\begin{split}
&  \int_{G} \phi(gx_0) (\tau(g) u, v) d g \\
=& \int_{[g] \in G/G_0}  \int_{G_0}  \phi(g g_0 x_0) (\tau(g g_0)u, v) d g_0 d [g G_0]\\
=& \int_{y \in M} \phi(y)  \int_{G_0}(\tau(\gamma(y) g_0) u, v) d g_0 d y \\
=& \int_{M} \phi(y)  \int_{G_0} (\tau(g_0)u, \tau(\gamma(y)^{-1})v) d g_0 dy \\
=&  \int_{G_0} (\tau(g_0)u,  \int_{M} \overline{\phi(y)} \tau(\gamma(y)^{-1}) v dy) d g_0\\
=&  \int_{G_0} (\tau(g_0)u,  v_0) d g_0.
\end{split}
\end{equation}
Q.E.D.\\
\\
We can further utilize the right invariance of the Haar measure on $G$ by
 changing $x_0$ into an arbitrary $x \in M$. 
\begin{thm}~\label{orbit1} Under the same assumptions from Theorem ~\ref{orbit0},
suppose $\psi(x^{ })$ is an absolutely integrable function on $M$. Let
$$u_0=\int_{M} \overline{\psi(x^{ })} \tau(\gamma(x^{ })^{-1})u d x^{ }$$
Suppose the function
$$\phi(g^{-1} x^{ }) \overline{\psi(x^{ })}  (u, \tau(g) v)$$
is in $L^1(G \times M)$.
Then we have
$$\int_M \int_G \phi(g^{-1} x^{ }) \overline{\psi(x^{ })}  (u, \tau(g) v) d g d x^{ } =\int_{G_0} (\tau(g_0) u_0, v_0) d g_0 .$$
\end{thm}
Proof: First of all, since $\tau$ is unitary and $\psi(x^{ })$ is integrable, $u_0$ is well-defined. According to Fubini's theorem, we can interchange the order of integrations.
We obtain
\begin{equation}
\begin{split}
 & \int_M \int_G \phi(g^{-1} x^{ }) \overline{\psi(x^{ })}  (u, \tau(g) v) d g d x^{ } \\
 = & \int_{M} \int_{G} \phi(g x^{ }) \overline{\psi(x^{ })}  ( \tau(g)u, v) d g d x^{ }  \\
= & \int_{M} \overline{\psi(x^{ })} \int_{G} \phi(g x^{ }) (\tau(g)u, v) d g d x^{ } \\
=& \int_{M} \overline{\psi(x^{ })} (\int_{G} \phi(g \gamma(x^{ }) x_0) (\tau(g)u, v) d g) d x^{ } \\
=& \int_{M} \overline{\psi(x^{ })} (\int_{G} \phi(g x_0) (\tau(g \gamma(x^{ })^{-1})u, v) d g ) d x^{ } \, \, \, \mbox{by the right invariance of $dg$} \\
=& \int_{M} \overline{\psi(x^{ })} (\int_{G} \phi(g x_0) (\tau(g)\tau(\gamma(x^{ })^{-1})u, v) d g ) d x^{ } \\
=& \int_{M} \overline{\psi(x^{ })} \int_{G_0} (\tau(g_0)\tau( \gamma(x^{ })^{-1}) u, v_0)d g_0 d x \,\,\, \mbox{by Theorem ~\ref{orbit0}} \\
=& \int_{G_0} (\tau(g_0)(\int_{M} \overline{\psi(x^{ })} \tau(\gamma(x^{ })^{-1})u d x^{ }), v_0) d g_0 \\
=& \int_{G_0} (\tau(g_0) u_0, v_0) d g_0. \\
\end{split}
\end{equation}
Q.E.D. 
\subsection{Orbital Integral $I(\phi,u,\mc O_x)$}
Let $\mc O_x$ be a generic $G_1$-orbit in $X_{00}^0$.
Then $\mc O_x$ possesses an $G_1$-invariant
measure.
Let $\pi$ be a unitary representation in the semistable range of
$\theta(MG_1, MG_2)$. Let us recall some notations and facts 
from Part II.
\begin{enumerate}
\item $\xi$ is a central unitary character of $M^0 G_1$ and any element $g^0$ in $M^0 G_1$ can be expressed as a pair $(\xi, g)$ with $g$ in $G_1$. 
\item $\pi_0=\omega(M^{\p} G_1, M^{\p} G_2^{\p})^c \otimes \pi \otimes \overline{\xi}$ is a representation of $G_1$.
\end{enumerate}
 We fix a $K$-finite vector $u$ in $\pi \otimes \overline{\xi}$. Suppose
$$\phi=\sum_{i=1}^s \phi^0_i \otimes \phi^{\p}_i$$
with $\phi_i^0 \in \omega(M^{0} G_1, M^{0} G_2^{0})$ and
$\phi^{\p}_i \in \omega(M^{\p} G_1, M^{\p} G_2^{\p})$. 
Then we have
\begin{equation}~\label{or}
\begin{split}
& (\phi \otimes u, \phi \otimes u)_{\pi} \\
= & \int_{MG_1} (\omega(MG_1, MG_2)(\tilde g) \phi, \phi)(u, \pi(\tilde g) u) d \tilde g \\
= & \sum_{i,j}\int_{M^0 G_1} (\omega(M^0G_1, M^0 G_2^0)(g^0) \phi_i^0, \phi_j^0)
(\phi_j^{\p} \otimes u, (\omega(M^{\p} G_1, M^{\p} G_2^{\p} )^c \otimes \pi)(g^0) (\phi_i^{\p} \otimes u)) d g^0 \\
= & \sum_{i,j}\int_{M^0 G_1} ((\omega(M^0G_1, M^0 G_2^0)\otimes \overline{\xi})(g^0) \phi_i^0, \phi_j^0)
(\phi_j^{\p} \otimes u, (\omega(M^{\p} G_1, M^{\p} G_2^{\p} )^c \otimes \pi  \otimes \overline{\xi})(g^0) (\phi_i^{\p} \otimes u)) d g^0\\
=& 2 \sum_{i,j} \int_{G_1} \int_{X^0} \phi_i^0(g^{-1} x) \overline{\phi_j^0(x)}
d x ( \phi_j^{\p} \otimes u, \pi_0(g)(\phi_i^{\p} \otimes u)) d g \\
=& 2 \sum_{i,j} \int_{X^0} \int_{G_1} \phi_i^0(g^{-1} x) \overline{\phi_j^0(x)}
 (\phi_j^{\p} \otimes u , \pi_0(g)(\phi_i^{\p} \otimes u)) d g d x\\
\end{split}
\end{equation}
First of all, due to Theorem ~\ref{conver}, the above integral converges absolutely. 
Since $X_{00}^0$ is open and dense in $X^0$, 
$$ 2
\int_{X^0_{00}} \int_{G_1} \sum_{i,j}  \phi_i^0(g^{-1} x) \overline{\phi_j^0(x)}
 (\phi_j^{\p} \otimes u , \pi_0(g)(\phi_i^{\p} \otimes u)) d g d x$$
converges absolutetly. Due to Fubini's Theorem, for almost all the orbits $\mc O_{x}$ in $X_{00}^0$, the function
$$\phi_i^0(g^{-1} x) \overline{\phi_j^0(x)}  (\phi_j^{\p} \otimes u, \pi_0(g) (\phi_i^{\p} \otimes u)) \qquad \forall \ \ i,j \in [1,s]$$
is absolutely integrable on $\mc O_x \times G_1$. 
Secondly, since $\{\phi_j^0\}_{j=1}^s$ are rapidly decaying functions in the Schr\"odinger model of $\omega(M^0G_1, M^0G_2^0))$, $\{ \phi_j^0 \}_{j=1}^s$ are absolutely integrable on $X_{00}^0$.  Hence $\{\phi_j^0\}_{j=1}^s$ are absolutely integrable on almost every $G_1$ orbit $\mc O_x$. \\
\\
Take $M$ to be an $G_1$-orbit $\mc O_x$ such that
\begin{enumerate}
\item $\phi_j^0$ is absolutely integrable on $\mc O_x$ for every $j$;
\item The function 
$$\phi_i^0(g^{-1} x) \overline{\phi_j^0(x)}  (\phi_j^{\p} \otimes u, \pi_0(g) (\phi_i^{\p} \otimes u)) $$
is absolutely integrable on $\mc O_x \times G_1$ for every $ i,j \in [1,s]$.
\end{enumerate} 
Denote the orbital integral
$$ \sum_{i,j} \int_{\mc O_x} \int_{G_1} \phi_i^0(g^{-1} x) \overline{\phi_j^0(x)}
 (\phi_j^{\p} \otimes u ,\pi_0(g)(\phi_i^{\p} \otimes u)) d g d x$$
by $ I(\phi, u, \mc O_x) $.
Take $\tau$ to be $\pi_0$. 
Since $M$ can be identified with $G_1/{G_1}_x$, $G_1$ forms a fiber bundle over $M$. By local triviality, we choose a smooth section 
$\gamma: M \rightarrow G_1$
over an open dense subset of $M$. 
Then 
Theorem ~\ref{orbit1} implies
\begin{equation}
\begin{split}
 & I(\phi, u, \mc O_x) \\
 = & \sum_{i,j} \int_{\mc O_x} \int_{G_1} \phi_i^0(g^{-1} x) \overline{\phi_j^0(x)}
 (\phi_j^{\p} \otimes u ,\pi_0(g)(\phi_i^{\p} \otimes u)) d g d x \\
= & \sum_{i,j} \int_{{G_1}_x} (\pi_0(g_0) \int_{\mc O_x} \overline{\phi_j^0(y)} \pi_0(\gamma(y)^{-1})(\phi_j^{\p} \otimes u) d y, 
\int_{\mc O_x} \overline{\phi_i^0(y)} \pi_0(\gamma(y)^{-1})(\phi_i^{\p} \otimes u) d y) d g_0 \\
=& \int_{{G_1}_x} (\pi_0(g_0) u_0,u_0) d g_0. \\
\end{split}
\end{equation} 
Here
 $$u_0=\int_{\mc O_x} \sum_{i} \overline{\phi_i^0(y)} \pi_0(\gamma(y)^{-1})(\phi_i^{\p} \otimes u) d y.$$

\subsection{Compactly Supported Continuous Functions}
The theorems we have so far proved hold for compactly supported
 continuous (not necessarily smooth) functions $\phi_i^0, \psi_i^0$ as well. In fact, any compactly
supported continuous function on $X^0$ can be dominated by a multiple of the
Gaussian function $\mu(x)$ on $X^0$. Therefore, the function
$$ |\phi^0_i(g^{-1} x) \overline{\psi_j^0(x)}  (\psi_j^{\p} \otimes u, \pi_0(g) (\phi_i^{\p} \otimes v) )|  $$
 is always in $L^1(G_1 \times X^0)$. The rest of the argument from Part III
goes through. Again, we obtain
\begin{thm}~\label{basic}
Let $\pi$ be a unitary representation in the semistable range of
$\theta(MG_1, MG_2).$
Let $u$ be a $K$-finite vector in $\pi \otimes \overline{\xi}$.
 Let $\phi_i^0$ be compactly supported continuous functions on $X^0$ and
$\phi^{\p}_i \in \omega(M^{\p} G_1, M^{\p} G_2^{\p})$. Write 
$$\phi=\sum_{i=1}^s \phi^0_i \otimes \phi^{\p}_i.$$
Then the integral
$$(\phi \otimes u, \phi \otimes u)_{\pi}=2 \sum_{i,j} \int_{X^0_{00}} \int_{G_1} \phi_i^0(g^{-1} x) \overline{\phi_j^0(x)}
(\phi_j^{\p} \otimes u , \pi_0(g)(\phi_i^{\p} \otimes u)) d g d x $$
 is absolutely convergent.
For almost every $G_1$-orbit $\mc O$ (except a subset of measure zero), 
$I(\phi, u, \mc O)$ converges absolutely. Fix such an orbit $\mc O_x$ and a base point $x$. Choose any smooth section 
$\gamma: \mc O_x \rightarrow G_1$
over an open dense subset of $\mc O_x$. Let
$$u_0=\int_{\mc O_x} \sum_{i} 
\overline{\phi_i^0(y)} \pi_0(\gamma(y)^{-1})(\phi_i^{\p} \otimes u) d y.$$
Then
$$I(\phi, u, \mc O_x) = \int_{{G_1}_x}
 (\pi_0(g) u_0,u_0) d g. $$
\end{thm}

\section{Part IV: Positivity and Unitarity}
\begin{lem} Suppose $\pi$ is a unitary representation in $ \mc R(MG_1, \omega(MG_1, MG_2))$.
Suppose for every $\phi \in \omega(MG_1, MG_2)$ and a fixed nonzero 
$u \in \pi$
$$(\phi \otimes u, \phi \otimes u)_{\pi} \geq 0.$$
Then $(,)_{\pi}$ is positive semidefinite. If $(,)_{\pi}$ does not vanish,
Then $\theta(MG_1, MG_2)(\pi)$ is unitary.
\end{lem}
A similar statement can be found in ~\cite{pr1}. \\
\\
Proof: If $(,)_{\pi}$ vanishes, the lemma holds automatically.
Suppose $(,)_{\pi}$ does not vanish. Let $\mc R_{\pi}$ be the radical of $(,)_{\pi}$.
The linear space 
$$(\mc P \otimes u)/(\mc R_{\pi} \cap (\mc P \otimes u))$$
must be nontrivial. Otherwise $\mc P \otimes u \subseteq \mc R_{\pi}$. 
Since $\mc R_{\pi}$ is a $(\f g_1, MK_1)$-module, by the  $(\f g_1, MK_1)$-action,
$$\mc P \otimes \pi^c  \subseteq \mc R_{\pi}.$$
This contradicts the nonvanishing of $(,)_{\pi}$. \\
\\
Observe that
$$(\mc P \otimes u)/(\mc R_{\pi} \cap (\mc P \otimes u))$$
is an admissible Harish-Chandra module of $MG_2$.
From Theorem 7.8 ~\cite{theta}, it must be irreducible and equivalent to
$\mc P \otimes \pi^c /\mc R_{\pi}$. Since
$$\int_{MG_1} (\phi, \omega(g)\phi) (\pi(g) u, u) d g \geq 0$$
for a fixed $u \in \pi$ and any $K$-finite $\phi$, $(,)_{\pi}|_{\mc P \otimes u}$ induces
an invariant positive definite form on $\theta(MG_1, MG_2)(\pi)$. Thus $\theta(MG_1, MG_2)(\pi)$
must be unitary.
Consequently, $(,)_{\pi}$ must be positive semidefinite. Q.E.D.

\subsection{Proof of the Main Theorem}
\begin{thm}~\label{main}
Let $\Xi(g)$ be Harish-Chandra's basic spherical function of $G_1$. Suppose
 \begin{enumerate} 
\item  $\pi$ is a unitary representation in the semistable range of
$\theta(MG_1, MG_2)$.
\item For any $x,y \in G_1$, the function $\Xi(xgy)$ is
integrable on ${G_1}_{\phi}$ for every
 generic $\phi \in Hom_D(V_1, X_2^0)$ (see Definition ~\ref{gen}). 
\item $\pi_0$ is weakly contained in $L^{2}(G_1)$. 
\end{enumerate}
Then $(,)_{\pi}$ is positive semidefinite. If $(,)_{\pi}$ does
not vanish, then $\theta(MG_1, MG_2)(\pi)$ is unitary.
\end{thm}
Roughly speaking, the second condition requires  ${G_1}_{\phi}$ be half the "size"
of $G_1$. The first condition is redundant assuming the second and the third conditions are true. The third conditions can be converted into a growth condition
on the matrix coefficients of $\pi$. \\
\\
Proof of the Theorem: Let $u$ be a fixed $K$-finite vector in $\pi \otimes \overline{\xi}$. Write
$$\mc S=\{\phi=\sum_{i=1}^s \phi^0_i \otimes \phi^{\p}_i \mid \phi_i^0 \in C^{\infty}_c(X^0), \phi_i^{\p} \in \omega(M^{\p} G_1, M^{\p} G_2^{\p}) \}.$$
Let $\phi \in \mc S$. 
Choose an arbitrary $G_1$-orbit $\mc O_x$ in $X_{00}^0$ such that $I(\phi, u, \mc O_x)$ converges absolutely.
There is a canonical fibre bundle
$${G_1}_{x} \rightarrow G_1 \rightarrow \mc O_x.$$
Fix a smooth section 
$\gamma: \mc O_x \rightarrow G_1$
over an open dense subset of $\mc O_x$ such that
the closure of $\gamma(supp(\phi_i^0))$ is compact for every $i$. Let
 $$u_0=\int_{\mc O_x} \sum_{i} \overline{\phi_i^0(y)} \pi_0(\gamma(y)^{-1})(\phi_i^{\p} \otimes u) d y.$$
From Theorem ~\ref{pos}, we have
$$\int_{{G_1}_x} (\pi_0(g) u_0,u_0) d g \geq 0 .$$
Combined with Theorem ~\ref{basic}, we obtain
$$I(\phi,u,\mc O_x) \geq 0;$$
$$(\phi \otimes u, \phi \otimes u)_{\pi} =\int_{\mc O \in G_1 \backslash X_{00}^0} I(\phi, u, \mc O) d [\mc O]  \geq 0.$$
We have thus proved that the Hermitian form $(,)_{\pi}$ restricted to $\mc S \otimes u$ is
positive semidefinite, i.e.,
$$\int_{MG_1}(\omega(MG_1, MG_2)(\tilde g) \phi, \phi)(u, \pi(\tilde g) u) d \tilde g \geq 0 \\$$
for every $\phi \in \mc S$. \\
\\
For an arbitrary $K$-finite vector $f$ in $\omega(MG_1, MG_2)$,  write
$$f=\sum_{k=1}^{s} f_k^0 (x) \otimes f_k^{\p} \qquad
(f_k^0 \in \omega(M^0 G_1, M^0 G_2^0), f_k^{\p} \in \omega(M^{\p} G_1, M^{\p} G_2^{\p})).$$
For each $k$, choose a sequence
 $\psi_k^{(j)}(x) \in C_c^{\infty}(X^0)$ such that
$$|\psi_k^{(j)}(x)| \leq |f(x)|$$
$$\psi_k^{(j)}(x) \rightarrow f(x).$$
Let $\psi^{(j)}=\sum_{k=1}^{s} \psi_k^{(j)} \otimes f_k^{\p}$. Apparently,
$\psi^{(j)} \in \mc S$ and
$$(\omega(MG_1, MG_2)(\tilde g) \psi^{(j)}, \psi^{(j)})(u, \pi(\tilde g) u) 
\rightarrow (\omega(MG_1, MG_2)(\tilde g) f, f)(u, \pi(\tilde g) u) $$
pointwise. Furthermore,
$$|(\omega(MG_1, MG_2)(\tilde g) \psi_j, \psi_j)(u, \pi(\tilde g) u) |
\leq \sum_{k,i=1}^s |(\omega(MG_1, MG_2)(\tilde g) |f_k^{0}| \otimes f_k^{\p}, |f_i^{0}| \otimes f_i^{\p})(u, \pi(\tilde g) u) |.$$
By the definition of semistable range, the function
$$|((\omega(MG_1, MG_2)(\tilde g)| f_k^0| \otimes f_k^{\p}, |f_i^0| \otimes f_i^{\p})(u, \pi(\tilde g) u) |$$
is absolutely integrable on $MG_1$ (see ~\cite{theta}). Hence, by dominated convergence theorem, 
$$(f \otimes u, f \otimes u)_{\pi}=\lim_{j \rightarrow \infty} (\psi^{(j)} \otimes u, \psi^{(j)} \otimes u)_{\pi} \geq 0$$
Therefore, the form $(,)_{\pi}$ is positive semidefinite. If $(,)_{\pi}$ does not vanish, then $(,)_{\pi}$ considered as a form on
$$\theta(MG_1, MG_2)(\pi)$$
is positive definite ( see ~\cite{theta}). 
We conclude that $\theta(MG_1, MG_2)(\pi)$ is unitary. Q.E.D. \\
\\
For $(G_1,G_2)$ in the stable range, the generic isotropic group
${G_1}_{\phi}$ will be trivial. In this case, if $\pi$ is an irreducible unitary representation of $MG_1$, then $(,)_{\pi}$ is positive semidefinite and nonvanishing. This result is due to Li (~\cite{li2}).
\subsection{$G_1=Sp_{2n}(\mb R)$}
Take $G_1=Sp_{2n}(\mb R)$ as an example. We can make our theorem more precise.
First let me define a partial order $\preceq$ in $\mb R^n$.
We say that $a \preceq b$ if and only if
$$\sum_{j=1}^k a_j \leq \sum_{j=1}^k b_j$$
for all $k$.

\begin{cor}~\label{co1}
Suppose $n < p \leq q$. Let $\pi$ be an irreducible unitary
representation of $MSp_{2n}(\mb R)$. Suppose for every leading exponent (see Chapter 8.8 in ~\cite{knapp}) $v$ of $\pi$ we have
$$\Re(v) -(\frac{p+q}{2}-n-1) \preceq -\rho(Sp_{2n}(\mb R)).$$
Then $(,)_{\pi}$ is positive semidefinite. In addition, if $(,)_{\pi}$ is nonvanishing, then
 $$\theta(MG_1, MG_2)(\pi)$$
 is unitary. 
\end{cor}
Proof: Take $V_1= \mb R^{2n}$ and $X_2^0=\mb R^{n+1}$. Then
$V_2^{\p}$ is a linear space equipped with a nondegenerate symmetric form
of signature $(p-n-1,q-n-1)$. We verify the conditions in Theorem ~\ref{main}.
\begin{itemize}
\item For $x \in \Hom(V_1, X_2^0)$, the generic isotropic group ${G_1}_{x}$
is just $Sp_{n-1}(\mb R)$ for $n$ odd. For $n$ even,
the generic ${G_1}_{x}$ can be identified with $Sp_{n-2}(\mb R) \times N$ where $N \cong \mb R^{n}$. One can easily check that $\Xi(g)$ for $Sp_{2n}(\mb R)$
is integrable on ${G_1}_{x}$.
\item Since
$$\Re(v) -(\frac{p+q}{2}-n-1) \preceq -\rho(Sp_{2n}(\mb R)),$$
$\pi_0=\omega(M^{\p} G_1, M^{\p} G_2^{\p})^c \otimes \pi \otimes \overline{\xi}$
has almost square integrable matrix coefficients. According to Theorem 1 of
~\cite{chh}, $\pi_0$ is weakly contained in $L^2(G_1)$.
\item The first condition is really redundant.  By Theorem 3.2 ~\cite{li2},
matrix coefficients of $\omega(MO(n+1,n+1), MSp_{2n}(\mb R))$ are in
$L^{2-\delta}(MSp_{2n}(\mb R))$ for small $\delta >0$. Since $\pi_0$ is almost square integrable, the matrix coefficients of $\omega(MO(p,q), MSp_{2n}(\mb R))
\otimes \pi$ are in $L^{1-\delta_0}(MG_2)$ for  small $\delta_0 >0$. Thus
$\pi$ must be in the semistable range of $\theta(MG_1, MG_2)$.
\end{itemize}
We conclude that $(,)_{\pi}$ is positive semidefinite. Q.E.D. 
\subsection{$G_1=O(p,q)$}
Similarly, we obtain
\begin{cor}
Suppose $p+q \leq 2n+1$. Let $\pi$ be an irreducible unitary representation of $MO(p,q)$. Suppose for every leading exponent $v$ of $\pi$ we have
$$\Re(v) -(n-\frac{p+q}{2}) \preceq - \rho(O(p,q)).$$
Then $(,)_{\pi}$ is positive semidefinite. In addition, if $(,)_{\pi}$ is nonvanishing, then
 $$\theta(MG_1, MG_2)(\pi)$$
 is unitary. 
\end{cor}
For $p+q$ odd, the growth condition concerning the leading exponent $v$ can be strengthened to allow
$$\Re(v) -(n-\frac{p+q-1}{2}) \preceq - \rho(O(p,q)).$$
The proof is omitted.

\end{document}